\documentclass[10pt]{article}

\usepackage{latexsym}
\usepackage{amssymb} 
\usepackage{epsfig}
\usepackage{amsmath}
\usepackage[mathscr]{eucal}

\newtheorem{Lemma}{Lemma}
\newtheorem{Proposition}[Lemma]{Proposition}
\newtheorem{Theorem}[Lemma]{Theorem}

\newtheorem{Corollary}[Lemma]{Corollary}

\newcommand{\cd}{\ \stackrel{d}{\longrightarrow} \ }

\newcommand{\FF}{\mbox{${\mathcal F}$}}

\newcommand{\eps}{\varepsilon}
\newcommand{\bE}{{\bf E}}
\newcommand{\bP}{{\bf P}}

\newcommand{\bR}{{\bf R}}
\newcommand{\bB}{{\bf B}}

\newcommand{\Pbb}{\mbox{${\mathbb P}$}}

\newcommand{\sfrac}[2]{{\textstyle\frac{#1}{#2}}}

\newcommand{\qed}{\hfill{\ \ \rule{2mm}{2mm}} \vspace{0.2in}}
\newcommand{\proof}{\noindent \emph{Proof:}\ }

\newcommand{\Zbold}{{\mathbb{Z}}}

\title{On the One Dimensional Critical ``Learning from Neighbours" Model}
\author{
{\bf Antar Bandyopadhyay},
\footnote{E-Mail: {\tt antar@isid.ac.in}} \\
\and
{\bf Rahul Roy},
\footnote{E-Mail: {\tt rahul@isid.ac.in}} \\
\and
{\bf Anish Sarkar}
\footnote{E-mail: {\tt anish@isid.ac.in}} \\ 
\ \\
Theoretical Statistics and Mathematics Unit \\
Indian Statistical Institute, New Delhi\\
7 S. J. S. Sansanwal Marg \\
New Delhi 110016 \\
INDIA
}

\begin{document}

\maketitle

\begin{abstract}

We consider a model of a discrete time ``interacting particle system" on the 
integer line where infinitely many changes are allowed at each instance of time. We describe the model
using chameleons of two different colours, {\it viz}\/., red ($R$) and blue ($B$). 
At each instance of time each chameleon performs an
independent but identical coin toss experiment with probability $\alpha$ to decide whether to change its colour or not. 
If the coin lands head then the creature retains its colour
(this is to be interpreted as a ``success"), otherwise it observes the colours and coin tosses of its two nearest
neighbours and changes its colour only if, among its neighbors and including itself,   the proportion of successes of the
other colour is larger than the  proportion of successes of its own colour. This produces a Markov chain with
infinite state space $\left\{R, B\right\}^{\Zbold}$. 
This model was studied by Chatterjee and Xu \cite{ChXu04} in the context of diffusion of technologies in a set-up of
myopic, memoryless agents. In their work they assume different success probabilities of coin tosses according to the colour of the chameleon. 
In this work we consider the ``critical" case where the success probability, $\alpha$,  is the same irrespective of the colour of the chameleon. We show that starting from any initial translation invariant distribution of
colours the Markov chain converges to a limit of a single colour, i.e.,  even at the critical case there is no ``coexistence" of the two colours at the limit. As a
corollary we also characterize the set of all translation invariant stationary laws of this
Markov chain. Moreover we show that starting with an i.i.d. colour distribution with
density $p \in [0,1]$ of one colour (say red), the limiting distribution is all red with
probability $\pi\left(\alpha, p\right)$ which is continuous in $p$ and for $p$ ``small" $\pi(p) >
p$. The last result can be interpreted as the model favours the
``underdog".
 
\vspace{0.1in}

\vspace{0.1in}

\end{abstract}

\section{Introduction and Main Results}
\label{Sec:Intro}

\subsection{Background and Motivation}
\label{Subsec:Background}
Chatterjee and Xu \cite{ChXu04} introduced a model of particle systems consisting of a countable number of particles of two types, each particle situated on integer points of the integer line. The type of a particle evolves with time depending on the behaviour of the neighbouring particles. This model, as Chatterjee and Xu explain is 
\begin{quote}
\noindent
\emph{
``... a problem of diffusion of technology, where one technology is better than the other and 
agents imitate better technologies among their neighbours.''
}
\end{quote}

The model above is part of a large class of models studied by economists over the last decade on `social learning'. Ellison and Fudenberg \cite{Ell93} introduced the notion of social learning -- they studied how  the speed of learning and the ultimate determination of market equilibrium is affected by
 social networks and other institutions governing communication between market participants. Bala and Goyal \cite{BG2000} studied a model where 
 \emph{``... individuals periodically make decisions concerning the continuation of existing information links and the formation of new information links, with their cohorts ... (based on) ...  the costs of forming and maintaining links against the potential rewards from doing so.''}. They studied the long run behaviour of this process. Much of the work on this was inspired by an earlier paper Bala and Goyal \cite{BG1998} where the learning was from neighbours and they showed  \emph{`` ... local learning ensures that all agents obtain the same payoffs in the long run.''} Banerjee and Fudenberg \cite{BF2004} also obtained similar results of a single `long-run outcome' when the decision making of an individual is based on a larger group of cohorts.

Here we consider the model studied by  Chatterjee and Xu \cite{ChXu04}. Instead of particles or technologies we describe the model with chameleons which can change 
their colours. Let $G = \left(V, E\right)$ be an infinite connected graph which is 
\emph{locally finite}, i.e., $\mbox{deg}_G\left(v\right) < \infty$ for any vertex $v \in V$. 
Suppose at every $v \in V$ and at any instance of time $t$, there is a chameleon $\xi_v(t)$, 
which is either red (R) or blue (B) in colour. In accordance with its colour it 
also has either a red coin $C_R(v,t)$ or a blue coin $C_B(v,t)$. 
The red coin has a probability $p_R$ of success $(1)$ and a probability $1-p_R$ of failure $(0)$, while the blue coin has a probability $p_B$ of success $(1)$ and a probability $1-p_B$ of failure $(0)$. The outcome of a coin of a chameleon is independent of the outcomes of the coins as well as the colours of the other chameleons. The evolution is governed by the rule described below which is referred as \emph{Rule-I} in \cite{ChXu04}.

Fix $t \geq 0$ and $v \in V$, let 
$N_v := \left\{ u \,\Big\vert\, <v,u> \in E \,\right\} \cup \left\{ v \right\}$ 
be the set of neighbours of the vertex $v$ including itself.
\begin{itemize}
\item If $C_{\xi_v(t)}(v,t) = 1$ then $\xi_v(t+1) = \xi_v(t)$, in other words, if the coin toss of the
      chameleon at $v$ at time $t$ results in a success then it retains its colour. 
\item If $C_{\xi_v(t)}(v,t) = 0$ then it retains its colour if the \emph{proportion of successes
      of the coin tosses of the chameleons of its colour in $N_v$ is larger or equal to the proportion of
      successes of the coin tosses of the chameleons of the other colour in $N_v$}. Otherwise it switches to
      the other colour. 
\end{itemize}


Formally, we have a configuration $\xi(t) \in \{R,B\}^V$ for every $t \geq 0$ 
and two independent collections 
$\{C_R(v,t): v \in V, t \geq 0\}$ and $\{C_B(v,t): v \in V, t \geq 0\}$ of
i.i.d. $\mbox{Bernoulli}\left(p_R\right)$ and i.i.d. $\mbox{Bernoulli}\left(p_B\right)$
random variables. Let $\left(\Omega, \FF, \Pbb\right)$ be the probability space
where all these random variables are defined. 

The process $\xi_0$ starts with some initial distribution $\Pbb_0$ on $\{R,B\}^V$ and the evolution is governed by the rule above. Let $\Pbb_t$ be the distribution of $\xi_t$ at time $t$.
In this work we are interested in finding the possible limiting distributions $\pi$ for 
$\left(\Pbb_t\right)_{t \geq 0}$.
From the definition it follows that $\{\xi_t : t \geq 0\}$ is a Markov chain with state space $\{R,B\}^V$; 
thus any limiting distribution $\pi$, if it exists, is a stationary distribution of this Markov chain.
We also observe that there are two absorbing states for this Markov chain, namely, 
the configuration of \emph{all reds} and the configuration of \emph{all blues}. 
Let $\delta_{\bR}$ denote the degenerate 
measure on $ \{R,B\}^V$ which assigns mass $1$ to the configuration of \emph{all reds}, and similarly
$\delta_{\bB}$ denote the measure on $\{R,B\}^V$ which assigns mass $1$ to the configuration 
of \emph{all blues}. Chatterjee and Xu \cite{ChXu04} studied this model for the one dimensional 
integer line $\Zbold$ with nearest neighbor links. They showed that 
when $\left(\xi_i(0)\right)_{i \in \Zbold}$ are i.i.d. with $\Pbb_0\left(\xi_0(0) = R\right) = p$ and
$\Pbb_0\left(\xi_0(0) = B\right) = 1 - p$, for some
$p \in \left(0,1\right)$ and $p_R > p_B$
\begin{equation}
\Pbb_t \mbox{ converges weakly to } \delta_{\bR} \mbox{ as } t \to \infty .
\label{Equ:Convergence-to-All-Red}
\end{equation}

In this work we first present a simpler proof of the above result. 
However our main interest is the study of the model when $p_R = p_B$, that is, 
when the success/failure of a coin is ``colour-blind''. 
We call this the ``critical case''.
The following subsection provides our main results.

\subsection{Main Results}
\label{Subsec:Results}

We first state the result of Chatterjee and Xu \cite{ChXu04} for which we provide a different proof in Section \ref{Sec:Red-over-Blue}. 

\begin{Theorem}
\label{Thm:ChXu} 
Let $G := \Zbold$ be the one dimensional integer line with nearest neighbour links and
suppose $\{\xi_i(0) :i \in \mathbb Z\}$ are i.i.d. 
with $\Pbb_0\left(\xi_0(0) = R\right) = p$ and
$\Pbb_0\left(\xi_0(0) = B\right) = 1 - p$, for some
$p \in \left(0,1\right)$ 
If $p_R > p_B$, then 
\begin{equation}
\Pbb_t \mbox{ converges weakly to } \delta_{\bR} \mbox{ as } t \rightarrow \infty .
\end{equation}
\end{Theorem}

Our main result is for  the ``critical'' case when $p_R = p_B$. For this we have the following result.  
\begin{Theorem}
\label{Thm:Main}
Let $G := \Zbold$ with nearest neighbour links and
suppose $\{\xi_i(0) :i \in \mathbb Z\}$ are i.i.d. 
with $\Pbb_0\left(\xi_0(0) = R\right) = p$ and
$\Pbb_0\left(\xi_0(0) = B\right) = 1 - p$, for some
$p \in \left(0,1\right)$
Assume $p_R = p_B = \alpha \in \left(0,1\right)$. Then, as $t \to \infty$
\begin{equation}
\Pbb_t \mbox{ converges weakly to } \pi\left(\alpha, p\right) \delta_{\bR} + \left(1-\pi\left(\alpha,p\right)\right) \delta_{\bB}  ,
\label{Equ:Convergence-to-All-Same-Color}
\end{equation}
where $\pi\left(\alpha, p\right) \in \left[0,1\right]$ satisfies the following properties
\begin{itemize}
\item[(i)]    For every fixed $\alpha \in \left(0,1\right)$ the function
              $p \mapsto \pi\left(\alpha,p\right)$ is continuous on $\left[0,1\right]$.  
\item[(ii)]   For any $\alpha \mbox{\ and\ } p \in \left(0,1\right)$, 
              $\pi\left(\alpha, p\right) = 1 - \pi\left(\alpha, 1-p\right)$. 
              Thus in particular
              $\pi\left(\alpha, \sfrac{1}{2}\right) = \sfrac{1}{2}$ for all 
              $\alpha \in  	\left(0,1\right)$. 
\item[(iii)]  $p^2 < \pi\left(\alpha, p\right) < 2p -p^2$ for all 
              $0 < p < 1$ and $\alpha \in \left(0,1\right)$. 
\item[(iv)]   For every $\alpha \in \left(0,1\right)$ there exists $\eps \equiv \eps\left(\alpha\right) > 0$ 
              such that $\pi\left(\alpha, p\right) > p$ for all $0 < p < \eps$. 
\end{itemize}
\end{Theorem}

Theorem \ref{Thm:Main} basically says that under the evolution scheme described 
above if $p_R = p_B = \alpha$ then starting with i.i.d. colours on the
integer line the distribution of the colours will converge either to \emph{all red} or to \emph{all blue}
configuration. Thus ruling out the possibility of any \emph{coexistence} of both the colours at the limit.
Such a result is expected considering the one dimensionality of the graph $\Zbold$. This lack of coexistence on $\Zbold$
is akin to the situation in many statistical physics models, such as, percolation, Ising model, 
$q$-Potts model which do not admit \emph{phase transition} in one dimension \cite{Gri99, Geo88}.

\begin{figure*}
\begin{center}\leavevmode
\epsfxsize=4cm \epsfbox{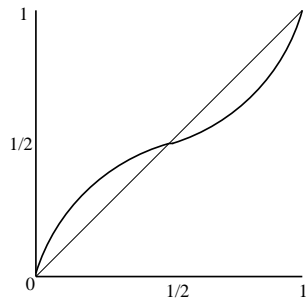}
\caption{Conjectured Graph of the function $p \mapsto \pi\left(\alpha, p\right)$.}
\end{center}
\label{Fig:pi-Curve}
\end{figure*}

It is interesting to note that 
$\pi\left(\alpha, p \right) > p$ in a neighbourhood of $0$, which can
be interpreted as follows:
\begin{quote}
The model gives an ``\emph{advantage to the underdog}'', in the sense that
for a fixed $\alpha \in \left(0,1\right)$ if $p$ is ``small'' then there is still a possibility that 
the (underdog) red chameleons will survive at the end. 
\end{quote}
We believe that this phenomenon is true for any $0 < p < 1$ with the caveat that the colour of the underdog is different according as $p$ is smaller or greater than $\sfrac{1}{2}$. We conjecture that the graph of the function
$p \mapsto \pi\left(\alpha, p\right)$ is as in Figure \ref{Fig:pi-Curve}
for every fixed $\alpha \in \left(0,1\right)$.


\subsection{Outline}
\label{Subsec:Outline}
The rest of the paper is divided as follows. In Section \ref{Sec:Red-over-Blue} we prove Theorem \ref{Thm:ChXu}. 
In Section \ref{Sec:One-Directional-Model} we consider a toy model on the half-line 
$\Zbold_+$ where a chameleon decides to change its colour according to its own toss and the colour and outcome of the toss of its neighbour to the right. 
This model is simpler to analyze and its usefulness is in providing an illustration of our method. In Section 
\ref{Sec:Proof-of-Main-Theorem} we prove Theorem \ref{Thm:Main}. 
A generalization of Theorem \ref{Thm:Main} is provided in
Section \ref{Sec:Translation-Invariant}, where we also characterize the translation invariant
stationary measures for this Markov chain. Section \ref{Sec:Technical} provides some auxiliary
technical results which we use in various derivations. 
We end with some discussion on the \emph{coexistence} of two 
colours in Section \ref{Sec:Coexistence}.

\section{Red is More Successful than Blue}
\label{Sec:Red-over-Blue}
In this section we provide a simple proof of the Theorem \ref{Thm:ChXu}. 
We begin by placing only blue chameleons at each point of the negative half line and red chameleons on the non-negative half line, and we take this as our initial configuration, that is,
\begin{equation}
\xi_i(0) = \begin{cases}
            B \mbox{ for } i \leq -1\\
R \mbox{ for } i \geq 0 .
           \end{cases}
\label{Equ:Blue-then-Red-Starting-Config}
\end{equation}

It is easy to see that, starting with $\xi(0)$ as given above there is always a sharp
interface left of which all chameleons are blue and right of which all are red. Moreover 
if we write $X_t$ as the position of the left most red chameleon at time $t \geq 0$ then
$\left(X_t\right)_{t \geq 0}$  performs a symmetric random walk starting from the
origin with i.i.d. increments, each taking values  $-1$, $0 $ and $1$ with probabilities
$p_R\left(1-p_B\right)$, $p_R p_B + \left(1-p_R\right) \left(1 - p_B\right)$ and
$p_B\left(1- p_R\right)$ respectively. This is because, for any $t \geq 0$ we have
\begin{equation}
X_{t+1} - X_t = \left\{\begin{array}{rll} 
                -1 & \mbox{if} & C_B\left(X_t-1,t\right) = 0 \mbox{\ and\ } C_R\left(X_t,t\right) = 1 \\
                 0 & \mbox{if} & C_B\left(X_t-1,t\right) = C_R\left(X_t,t\right) = 1 \mbox{\ or\ } \\
                   &           & C_B\left(X_t-1,t\right) = C_R\left(X_t,t\right) = 0\\
                +1 & \mbox{if} & C_B\left(X_t-1,t\right) = 1 \mbox{\ and\ } C_R\left(X_t,t\right) = 0
                \end{array}\right.
\label{Equ:Increment}
\end{equation}
Now it is easy to check that if $p_R > p_B$, then $\left(X_t\right)_{t \geq 0}$ is a transient
random walk with a strictly negative drift. In other words it proves that starting with the configuration
given in (\ref{Equ:Blue-then-Red-Starting-Config})
\[
\Pbb_t \mbox{ converges weakly to }\delta_{\bR} \mbox{ as } t \to \infty .
\]

\subsection{Proof of Theorem \ref{Thm:ChXu}}
\label{Subsec:Proof-of-Theorem-ChXu}
To complete the proof of Theorem \ref{Thm:ChXu},
fix an $\eta > 0$ and let $M$ be such that the random walk $X_t$ satisfies
\[
\Pbb\left(X_t \leq M \mbox{\ for all\ } t \geq 0\right) > 1- \eta  .
\]
If the initial distribution of chameleons at time $0$ is such that, for some $j \in \mathbb Z$
\[
\xi_i(0) = \begin{cases}
           R \mbox{ for } j \leq i \leq j+ 2M+1\\
           B \mbox{ otherwise} 
           \end{cases}
\]
then, 
\begin{equation}
\Pbb\left(X_t^- \rightarrow - \infty \mbox{\ and\ } X_t^+ \rightarrow \infty \right) 
\geq \left(1-\eta\right)^2
\end{equation}
where we write $X_t^-$ as the \emph{left interface} and $X_t^+$ as the \emph{right interface} at time 
$t \geq 0$.

Further an easy coupling argument shows that
the above situation of a stretch of $2M+2$ red chameleons flanked by only blue chameleons on either sides
is ``worse'' than the case when the two ends instead of being all blue is actually a mixture of red and blue chameleons. More precisely, suppose the starting configuration $\xi_0$ is such that 
there exists a location $j \in \Zbold$ with $\xi_i(0)=R$ for all $j - M \leq i \leq j+M+1$, then
\begin{equation}
\Pbb\left(X_t^{-,j} \rightarrow - \infty \mbox{\ and\ } X_t^{+,j} \rightarrow \infty \right) 
\geq \left(1-\eta\right)^2
\label{Equ:2-Interfaces}
\end{equation}
where we write $X_t^{-,j}$  and $X_t^{+,j}$ are the positions of the leftmost and rightmost red chameleons at time $t \geq 0$ of the  of the lounge (possibly enlarged) of chameleons which started as the stretch of length $2M+2$.

Let 
$J := \inf\left\{ j \geq 0 \,\Big\vert\, \xi_i(0) = R \mbox{\ for all\ } j - M \leq i \leq j+M+1 \,\right\}$.
When the initial distribution of the chameleons is i.i.d. then $\Pbb\left(0 \leq J < \infty\right) = 1$. 
Thus, conditioning on $ X_t^{-,J} \leq J - M$ and $ X_t^{+,J} \geq J+M+1$, the ensuing conditional independence of $ X_t^{-,J}$ and $ X_t^{+,J}$, yields
\begin{equation}
\Pbb\left(X_t^{-,J} \rightarrow - \infty \mbox{\ and\ } X_t^{+,J} \rightarrow \infty \right) 
\geq \left(1-\eta\right)^2 .
\label{Equ:2-Interfaces-IID}
\end{equation}
Now for any $k \geq 1$
\begin{eqnarray*}
\lefteqn{\liminf_{t \rightarrow \infty} \Pbb_t\left(\xi_i(t) = R \mbox{\ for all\ } -k \leq i \leq k \right)} \\
 & \geq & \liminf_{t \rightarrow \infty}
          \Pbb\left(X_t^{-,J} < -k \mbox{\ and\ } X_t^{+,J} > k \right) \\
 & \geq & \Pbb\left(X_t^{-,J} < -k \mbox{\ and\ } X_t^{+,J} > k  \mbox{\ eventually\ } \right) \\
 & \geq & \Pbb\left(X_t^{-,J} \rightarrow - \infty \mbox{\ and\ } X_t^{+,J} \rightarrow \infty \right) \\
 & \geq & \left(1-\eta\right)^2  ,
\end{eqnarray*}
where the last inequality follows from equation (\ref{Equ:2-Interfaces-IID}). 

Finally since $\eta > 0$ is arbitrary we conclude that $\Pbb_t$ converges weakly to $\delta_{\bR}$ as $ t \to \infty$.
$\qed$

\noindent
{\bf Remark:} We observe that the above argument holds for any starting configuration $\xi\left(0\right)$ such that
intervals of reds of arbitrary length can be found with probability one. This generalizes the 
Theorem \ref{Thm:ChXu}. 

\section{One Directional Neighbourhood Model}
\label{Sec:One-Directional-Model}
In this section we  study the simpler \emph{one directional neighbourhood model} model where the dynamics follows our rule but with $N_i := \left\{i, i+1\right\}$ for
$i \in \Zbold$.  The computations for this model are much simpler than the original \emph{two sided}
neighbourhood model and the method used here is illustrative of the method employed for the original two sided neighbourhood model. We now state the convergence result
for the one directional neighbourhood model. 

\begin{Theorem}
\label{Thm:One-Directional}
Let $\{\xi_i(0) : i \in \mathbb Z\}$ be i.i.d. random variables with 
$\Pbb(\xi_i(0) = R) = p = 1-\Pbb(\xi_i(0) = B)$. Then for the one directional neighbourhood model with
$p_R = p_B = \alpha \in \left(0,1\right)$ we have
\begin{equation}
\Pbb_t \cd p \delta_{\bR} + \left(1-p\right) \delta_{\bB} \mbox{\ \ as\ \ } t \rightarrow \infty  .
\label{Equ:One-Directional-Conv}
\end{equation}
\end{Theorem}
Before we prove this theorem we make the following observation which is very simple to prove but
plays important role in all our subsequent discussions.
\begin{Proposition}
\label{Prop:One-Directional-Translation-Inv}
Under the dynamics of the one directional neighborhood model, if
$\Pbb_0$ is a \emph{translation invariant} measure on $\left\{R, B\right\}^{\Zbold}$ then
$\Pbb_t$ is also translation invariant for every $t \geq 0$. 
\end{Proposition}
The proof of this Proposition follows from  the Markov chain dynamics of the model and we omit the details. 
It is worth remarking here that a similar result is true for the
two sided neighbourhood model.  

\subsection{Proof of Theorem \ref{Thm:One-Directional}}
\label{Subsec:Proof-of-One-Directional-Conv}
Before we embark on the proof of Theorem \ref{Thm:One-Directional} we present some notation.
Observe that, from the translation invariance of $\Pbb_t$ as given by Proposition \ref{Prop:One-Directional-Translation-Inv},  for every $t \geq 0$, $k \geq 1$, $i \in \Zbold$ and $\omega_j \in \{R, B\}$ 
$\Pbb_t\left(\xi_i(t)= \omega_1, \xi_{i+1}(t) = \omega_2, \ldots, \xi_{i+k-1} = \omega_k\right)$
does not depend on the location $i$; and thus with a slight ause of notation we write
$$
\Pbb_t\left(\omega_1, \omega_2, \ldots , \omega_k\right):= \Pbb_t\left(\xi_i(t)= \omega_1, \xi_{i+1}(t) = \omega_2, \ldots, \xi_{i+k-1} = \omega_k\right) .
$$
Also
$$
\Pbb_t \left(R\bB_kR\right):=\Pbb_t \left(RB\ldots BR\right)
$$
where there are $k$ many $B$'s in the expression on the right.

To prove this theorem we will use the technical result Theorem \ref{Thm:General-Conv} given in
Section \ref{Sec:Technical}.

Now fix $t \geq 0$. Observe
\begin{eqnarray}
\Pbb_{t+1} \left(R\right) 
 & = &   \Pbb_t\left(RR\right) + \left(\alpha + \left(1-\alpha\right)^2\right) \Pbb_t\left(RB\right)
       + \alpha \left(1-\alpha\right) \Pbb_t\left(BR\right) \nonumber \\
 & = &   \Pbb_t\left(RR\right) 
       + \left(\alpha + \left(1-\alpha\right)^2 + \alpha\left(1-\alpha\right)\right) \Pbb_t\left(RB\right)
         \nonumber \\
 & = &   \Pbb_t\left(RR\right) + \Pbb_t\left(RB\right) \nonumber \\
 & = &   \Pbb_t\left(R\right) \label{Equ:One-Directional-R}  .
\end{eqnarray}
The first equality follows from the dynamics rule. The second equality follows from the fact
$\Pbb_t\left(RB\right) = \Pbb_t\left(R\right) - \Pbb_t\left(RR\right) = \Pbb_t\left(BR\right)$, which is a consequence of the translation invariance of $\Pbb_t$.

Now for $t \geq 0$ using the rule of the dynamics we get
\begin{eqnarray}
\Pbb_{t+1}\left(RR\right) 
 & = & \Pbb_t\left(RRR\right) + \left(\alpha + \left(1-\alpha\right)^2\right) \Pbb_t\left(RRB\right) \nonumber \\
 &   & + \alpha \left(1-\alpha\right) 
         \left(\Pbb_t\left(RBR\right) + \Pbb_t\left(BRR\right) + \Pbb_t\left(BRB\right) \right) .
         \label{Equ:One-Directional-RR}
\end{eqnarray}
On the other hand by translation invariance of $\Pbb_t$ we have
\begin{equation}
\Pbb_t\left(RR\right) = \Pbb_t\left(RRR\right) + \Pbb_t\left(RRB\right)  .
\label{Equ:One-Directional-RR-Tran}
\end{equation}
Subtracting equation (\ref{Equ:One-Directional-RR-Tran}) from the
equation (\ref{Equ:One-Directional-RR}) we get
\begin{equation}
\Pbb_{t+1}\left(RR\right) - \Pbb_t\left(RR\right) 
= \alpha \left(1-\alpha\right) \left(\Pbb_t\left(RBR\right) + \Pbb_t\left(BRB\right) \right)  ,
\label{Equ:One-Directional-RR-RR}
\end{equation}
Here we use the fact that
$\Pbb_t\left(BRR\right) = \Pbb_t\left(RR\right) - \Pbb_t\left(RRR\right) = \Pbb_t\left(RRB\right)$. 
So we conclude that 
\begin{equation} 
\Pbb_t\left(RR\right) =  \Pbb_0\left(RR\right) 
                       + \alpha \left(1-\alpha\right) 
                         \sum_{n=0}^t \left(\Pbb_t\left(RBR\right) + \Pbb_t\left(BRB\right) \right)  .
\label{Equ:One-Directional-RR-Expression}
\end{equation}
Since the summands above are non-negative and since $0 \leq  \Pbb_t\left(RR\right) \leq 1$ we have  $\lim_{t \rightarrow \infty} \Pbb_t\left(RR\right)$ exists. In addition, we have
\[
\sum_{n=0}^{\infty} \left(\Pbb_t\left(RBR\right) + \Pbb_t\left(BRB\right) \right) < \infty  .
\]
So in particular
\begin{equation}
\lim_{t \rightarrow \infty} \Pbb_t\left(RBR\right) 
= 0 = 
\lim_{t \rightarrow \infty} \Pbb_t\left(BRB\right)  .
\label{Equ:One-Directional-RBR-BRB-Zero-Limit}
\end{equation}
Finally observe that for any $k \geq 0$ we have
\[
\alpha^{k+1} \left(1-\alpha\right) \Pbb_t\left(B\bR_kB\right) 
\leq 
\Pbb_{t+1}\left(B\bR_{k-1}B\right) \mbox{\ and} 
\]
\begin{equation}
\alpha^{k+1} \left(1-\alpha\right) \Pbb_t\left(R\bB_kR\right) 
\leq 
\Pbb_{t+1}\left(R\bB_{k-1}R\right)  .
\label{Equ:One-Directional-Mixed-Bound}
\end{equation}
Using
(\ref{Equ:One-Directional-RBR-BRB-Zero-Limit}) and (\ref{Equ:One-Directional-Mixed-Bound})
it follows by induction that 
\begin{equation}
\lim_{t \rightarrow \infty} \Pbb_t \left(R\bB_kR\right) 
= 0 = 
\lim_{t \rightarrow \infty} \Pbb_t \left(B\bR_kB\right)
\,\,\,\,
\forall \,\,\, k \geq 1  .
\label{Equ:One-Directional-Mixed-Zero-Limit}
\end{equation}
Theorem \ref{Thm:One-Directional} now follows from  Theorem \ref{Thm:General-Conv}. $\qed$

\subsection{Convergence from Translation Invariant Starting Distribution} 
\label{Subsec:One-Directional-General}
From the above proof of Theorem \ref{Thm:One-Directional}
we observe that nowhere have we used the fact that $\Pbb_0$ is an i.i.d. probability on
$\left\{R, B\right\}^{\Zbold}$, we just needed  $\Pbb_0$ to be a translation invariant measure. 
Thus the following generalization holds for this one directional neighborhood model. 
\begin{Theorem}
\label{Thm:One-Directional-General}
Let $\{\xi_i(0) : i \in \mathbb Z\}$ be random variables which are translation invariant and let
$\Pbb(\xi(0) = R) = p = 1-\Pbb(\xi_i(0) = B)$. Then for the one directional neighborhood model with
$p_R = p_B = \alpha \in \left(0,1\right)$ we have
\begin{equation}
\Pbb_t \cd p \delta_{\bR} + \left(1-p\right) \delta_{\bB} \mbox{\ \ as\ \ } t \rightarrow \infty  .
\end{equation}
\end{Theorem}
Moreover the following corollary is now immediate. 
\begin{Corollary}
\label{Cor:One-Directional-Characterization-Stationary-Measures}
For the one directional neighborhood model with $p_R = p_B = \alpha \in \left(0,1\right)$ 
the only translation invariant \emph{stationary} measures are of the form
\[
\kappa \delta_{\bR} + \left(1-\kappa\right) \delta_{\bB}
\]
for some $0 \leq \kappa \leq 1$. 
\end{Corollary}

\section{Proof of Theorem \ref{Thm:Main}}
\label{Sec:Proof-of-Main-Theorem}
In this section we will prove our main result, namely, Theorem \ref{Thm:Main}.
But before we proceed we note that as remarked in the previous section, 
the following result is also true for our original model.
\begin{Proposition}
\label{Prop:Translation-Inv}
Under the dynamics of our original model with $p_R=p_B$, if
$\Pbb_0$ is a \emph{translation invariant} measure on $\left\{R, B\right\}^{\Zbold}$ then
$\Pbb_t$ is also translation invariant for every $t \geq 0$. 
\end{Proposition}
Once again the proof is simple and hence we omit the details. 

As in Section \ref{Subsec:Proof-of-One-Directional-Conv}, Proposition
\ref{Prop:Translation-Inv} demonstrates the translation invariance of $\Pbb_t$ whenever $\Pbb_0$ is translation invariant.
The notation we use in this section are the same as set up in Section \ref{Subsec:Proof-of-One-Directional-Conv}.

\subsection{Proof of The Convergence (\ref{Equ:Convergence-to-All-Same-Color})}
As in the previous section, here too we use
the Theorem \ref{Thm:General-Conv} to prove the convergence 
(\ref{Equ:Convergence-to-All-Same-Color}). 
For that we begin by checking that
$\lim_{t \rightarrow \infty} \Pbb_t\left(RR\right)$ exists. In order to prove
this limit we use a similar technique as done in Section \ref{Sec:One-Directional-Model}. 
The dynamics of the two sided neighbourhood model bring in some additional intricacies. 

The following table presents some calculations which we use repeatedly. The  
column on the right is the probability of obtaining a configuration
$RR$ at locations $\left(i, i+1\right)$ of $\Zbold$ at time $(t+1)$ when the configuration at time $t$ at locations $\left(i-1, i, i+1, i+2\right)$ is given by the column on the left.
\begin{center}
\begin{tabular}{|c|cl|} \hline
Configuration at time $t$ & \qquad & Probability of getting a configuration $RR$ at time $t+1$ \\
\hline \hline
$RRRR$  &  \qquad & $1$ \\ \hline
$BRRR$  &  \qquad & $\alpha + \left(1-\alpha\right)^2$ \\ \hline
$RRRB$  &  \qquad & $\alpha + \left(1-\alpha\right)^2$ \\ \hline
$BRRB$  &  \qquad & $\alpha^2 + 2 \alpha \left(1-\alpha\right)^2 + \left(1-\alpha\right)^4$ \\ \hline
$RRBR$  &  \qquad & $\alpha \left(1-\alpha\right) \left(2 - \alpha\right)$ \\ \hline
$BRBR$  &  \qquad & $\alpha \left(1 - \alpha\right) \left(1 + \left(1-\alpha\right)^2\right)$ \\ \hline
$RRBB$  &  \qquad & $\alpha \left(1 - \alpha\right)$ \\ \hline
$BRBB$  &  \qquad & $\alpha \left(1 - \alpha\right)$ \\ \hline
$RBRR$  &  \qquad & $\alpha \left(1-\alpha\right) \left(2 - \alpha\right)$ \\ \hline
$BBRR$  &  \qquad & $\alpha \left(1 - \alpha\right)$ \\ \hline
$RBRB$  &  \qquad & $\alpha \left(1-\alpha\right) \left(2 - \alpha\right)$ \\ \hline
$BBRB$  &  \qquad & $\alpha \left(1 - \alpha\right)$ \\ \hline
$RBBR$  &  \qquad & $\alpha^2 \left(1-\alpha\right)^2$ \\ \hline
$BBBR$  &  \qquad & $0$ \\ \hline
$RBBB$  &  \qquad & $0$ \\ \hline
$BBBB$  &  \qquad & $0$ \\ \hline
\end{tabular}
\end{center}
Combining we get 
\begin{eqnarray}
\Pbb_{t+1}\left(RR\right) & = & 
\Pbb_t\left(RRRR\right) 
+ \left(\alpha + \left(1-\alpha\right)^2\right) 
  \left(\Pbb_t\left(BRRR\right) + \Pbb_t\left(RRRB\right)\right) \nonumber \\
 & & + \left( \alpha^2 + 2 \alpha \left(1-\alpha\right)^2 + \left(1-\alpha\right)^4 \right)
  \Pbb_t\left(BRRB\right) \nonumber \\
 & & + \alpha \left(1-\alpha\right) \left(2 - \alpha\right) 
  \left( \Pbb_t\left(RRBR\right) + Pbb_t\left(RBRR\right) \right) \nonumber \\
 & & + \alpha \left(1-\alpha\right) 
  \left(\Pbb_t\left(RRBB\right) + \Pbb_t\left(BRBB\right)
        \Pbb_t\left(BBRR\right) + \Pbb_t\left(BBRB\right)\right) \nonumber \\
& & + \alpha^2 \left(1-\alpha\right)^2 \Pbb_t\left(RBBR\right) \label{Equ:RR}
\end{eqnarray}
Also by translation invariance of $\Pbb_t$ it also follows that
\begin{equation}
\Pbb_t\left(RR\right) = 
\Pbb_t\left(RRRR\right) + \Pbb_t\left(RRRB\right) 
+ \Pbb_t\left(BRRR\right) + \Pbb_t\left(BRRB\right)
\label{Equ:RR-Tran}
\end{equation}
Now subtracting equation (\ref{Equ:RR-Tran}) from equation (\ref{Equ:RR}) we get
\begin{eqnarray}
\lefteqn{\Pbb_{t+1} \left(RR\right) - \Pbb_t\left(RR\right)} \nonumber \\
 & = & \left(\alpha + \left(1-\alpha\right)^2 - 1\right) 
  \left(\Pbb_t\left(BRRR\right) + \Pbb_t\left(RRRB\right)\right) \nonumber \\
 & & + \left( \alpha^2 + 2 \alpha \left(1-\alpha\right)^2 + \left(1-\alpha\right)^4 - 1\right) 
       \Pbb_t\left(BRRB\right) \nonumber \\
 & & + \alpha \left(1-\alpha\right) \left(2 - \alpha\right)
       \left( \Pbb_t\left(RRBR\right) + \Pbb_t\left(RBRR\right) \right) \nonumber \\
 & & + \alpha \left(1-\alpha\right) 
  \left(\Pbb_t\left(RRBB\right) + \Pbb_t\left(BRBB\right)
        \Pbb_t\left(BBRR\right) + \Pbb_t\left(BBRB\right)\right) \nonumber \\
 & & + \alpha^2 \left(1-\alpha\right)^2 \Pbb_t\left(RBBR\right) \nonumber \\
 & = & \alpha \left(1-\alpha\right) \left[
       -\left(\Pbb_t\left(RRRB\right) + \Pbb_t\left(BRRB\right) 
       + \Pbb_t\left(BRRR\right) + \Pbb_t\left(BRRB\right) \right) \right. \nonumber \\
 & & +\alpha \left(1 - \alpha\right) \left(\Pbb_t\left(BRRB\right) + \Pbb_t\left(RBBR\right)\right) \nonumber \\
 & & +\left(1-\alpha\right) \left(\Pbb_t\left(RRBR\right) + \Pbb_t\left(RBRR\right)\right) \nonumber \\
 & & +\left(1+\left(1-\alpha\right)^2\right) 
      \left(\Pbb_t\left(BRBR\right) + \Pbb_t\left(RBRB\right)\right) \nonumber \\
 & & +\left(\Pbb_t\left(BRBB\right) + \Pbb_t\left(RBRR\right)\right) \nonumber \\
 & & \left. +\left(\Pbb_t\left(RRBB\right) + \Pbb_t\left(RRBR\right) 
       + \Pbb_t\left(BBRR\right) + \Pbb_t\left(BBRB\right)\right) \right] \nonumber \\
& = & 
\alpha \left(1-\alpha\right) 
\left[
\alpha \left(1 - \alpha\right) \left(\Pbb_t\left(BRRB\right) + \Pbb_t\left(RBBR\right)\right)\right. \nonumber \\
& & +\left(1-\alpha\right) \left(\Pbb_t\left(RRBR\right) + \Pbb_t\left(RBRR\right)\right) \nonumber \\
& & +\left(1+\left(1-\alpha\right)^2\right) 
 \left(\Pbb_t\left(BRBR\right) + \Pbb_t\left(RBRB\right)\right) \nonumber \\
& & \left. +\left(\Pbb_t\left(BRBB\right) + \Pbb_t\left(BBRB\right)\right) \right] \label{Equ:RR-RR}
\end{eqnarray}

The last equality follows from the following:
\[ \Pbb_t\left(RRB\right) = \Pbb_t\left(RRBB\right) + \Pbb_t\left(RRBR\right)
                          = \Pbb_t\left(RRRB\right) + \Pbb_t\left(BRRB\right)  , \]
\[ \Pbb_t\left(BBR\right) = \Pbb_t\left(BBRR\right) + \Pbb_t\left(RBRR\right)
                          = \Pbb_t\left(BRRR\right) + \Pbb_t\left(BRRB\right)  .\]

Thus we have for any $t \geq 0$,
\begin{eqnarray}
\Pbb_{t+1}\left(RR\right) - \Pbb_0\left(RR\right) 
 & = &   \alpha \left(1-\alpha\right) \left[ \alpha \left(1 - \alpha\right) 
         \sum_{n=0}^t \left(\Pbb_n\left(BRRB\right) + \Pbb_n\left(RBBR\right)\right) \right. \nonumber \\
 &   & + \left(1-\alpha\right) \sum_{n=0}^t \left(\Pbb_n\left(RRBR\right) + \Pbb_n\left(RBRR\right)\right)
         \nonumber \\
 &   & + \left(1+\left(1-\alpha\right)^2\right) 
         \sum_{n=0}^t \left(\Pbb_n\left(BRBR\right) + \Pbb_n\left(RBRB\right)\right) \nonumber \\
 &   & + \left. \sum_{n=0}^t \left(\Pbb_n\left(BRBB\right) + \Pbb_n\left(BBRB\right)\right) \right]
       \label{Equ:RR-Expression}
\end{eqnarray}
Since all the terms and summands on the right of the above equality are non-negative, we have $\lim_{t \rightarrow \infty} \Pbb_t\left(RR\right)$ exists. Moreover we obtain that the sequence $\{\Pbb_t\left(\omega\right) : t \geq 0\}$ is
summable whenever 
\[ 
\omega \in \left\{ BRRB, RBBR, RRBR, RBRR, BRBR, RBRB, BRBB, RBRR \right\}  .
\] 
Using translation invariance we also have that the sequences
$\{\Pbb_t\left(RBR\right) : t \geq 0 \}$ and
$\{\Pbb_t\left(BRB\right) : t \geq 0 \} $ are summable.
In particular we conclude
\begin{equation}
\lim_{t \rightarrow \infty} \Pbb_t\left(RBR\right) 
= 0 = 
\lim_{t \rightarrow \infty} \Pbb_t\left(BRB\right)  .
\label{Equ:RBR-BRB-Zero-Limit}
\end{equation}

Now observe that for any $k \geq 0$ we have
\[
\alpha^{k+1} \left(1-\alpha\right) \Pbb_{t+1}\left(B\bR_kB\right) 
\leq 
\Pbb_t\left(B\bR_{k-1}B\right) \mbox{\ and} 
\]
\begin{equation}
\alpha^{k+1} \left(1-\alpha\right) \Pbb_{t+1}\left(R\bB_kR\right) 
\leq 
\Pbb_t\left(R\bB_{k-1}R\right)  .
\end{equation}

Finally, we consider the one dimensional marginal and observe
\begin{eqnarray}
\Pbb_{t+1}\left(R\right) 
 & = & \Pbb_t\left(RRR\right) + \left(\alpha + \left(1-\alpha\right)^2\right) 
                               \left(\Pbb_t\left(RRB\right) + \Pbb_t\left(BRR\right)\right) \nonumber \\
 &   & + \left(\alpha + \left(1  - \alpha\right)^3\right) \Pbb_t\left(BRB\right) \nonumber \\ 
 &   & + \alpha \left(1-\alpha\right) \left(\Pbb_t\left(BBR\right) + \Pbb_t\left(RBB\right)\right) \nonumber \\
 &   & + \left(1-\alpha\right) \left(1 - \left(1-\alpha\right)^2\right) \Pbb_t\left(RBR\right)
       \label{Equ:R}
\end{eqnarray}
Also from translation invariance of $\Pbb_t$ it follows that
\begin{equation}
\Pbb_t\left(R\right)
= 
\Pbb_t\left(RRR\right) + \Pbb_t\left(RRB\right) + \Pbb_t\left(BRR\right) + \Pbb_t\left(BRB\right)
\label{Equ:R-Tran}
\end{equation}
Subtracting equation (\ref{Equ:R-Tran}) from equation (\ref{Equ:R}) we have
\begin{equation}
\Pbb_{t+1}\left(R\right) - \Pbb_t\left(R\right) 
 = 
\alpha^2 \left(1-\alpha\right) \left(\Pbb_t\left(BRB\right) - \Pbb_t\left(RBR\right)\right)
\label{Equ:R-R}
\end{equation}
To derive this final expression we use the following identities which are easy consequences of
translation invariance of $\Pbb_t$.
\[
  \Pbb_t\left(BBR\right) - \Pbb_t\left(BRR\right) 
= \Pbb_t\left(BRB\right) - \Pbb_t\left(RBR\right)
= \Pbb_t\left(RBB\right) - \Pbb_t\left(RRB\right) 
\]
The summability of the sequences
$\{\Pbb_t\left(RBR\right): t \geq 0 \}$ and
$\{\Pbb_t\left(BRB\right): t \geq 0 \} $ yields, from equation \ref{Equ:R-R}, the existence  of  $\lim_{t \rightarrow} \Pbb_t\left(R\right)$.

Invoking Theorem \ref{Thm:General-Conv} we now complete the proof of the convergence (\ref{Equ:Convergence-to-All-Same-Color}).

\subsection{Proof of the Properties of $\pi\left(\alpha, p\right)$}
First, from the definition it follows that
$\pi\left(\alpha, p\right) = \lim_{t \rightarrow \infty} \Pbb_t\left(RR\right)$; thus using equation
(\ref{Equ:RR-Expression}) we get
\begin{eqnarray}
\pi\left(\alpha, p\right) 
 & = & p^2 + \alpha \left(1-\alpha\right) 
       \left[ \alpha \left(1 - \alpha\right) 
       \sum_{t=0}^{\infty} \left(\Pbb_t\left(BRRB\right) + \Pbb_t\left(RBBR\right)\right) \right. \nonumber \\
 &   & + \left(1-\alpha\right) 
         \sum_{t=0}^{\infty} \left(\Pbb_t\left(RRBR\right) + \Pbb_t\left(RBRR\right)\right) \nonumber \\
 &   & + \left(1+\left(1-\alpha\right)^2\right) 
         \sum_{t=0}^{\infty} \left(\Pbb_t\left(BRBR\right) + \Pbb_t\left(RBRB\right)\right) \nonumber \\
 &   & + \left. \sum_{t=0}^{\infty} \left(\Pbb_t\left(BRBB\right) + \Pbb_t\left(BBRB\right)\right) \right]   .
         \label{Equ:pi-Expression}
\end{eqnarray}
This immediately proves that $\pi\left(\alpha, p \right) > p^2$ for any $p \in \left(0,1\right)$. Moreover because 
the model is symmetric with respect to colour we have
\begin{equation}
\pi\left(\alpha, p\right) = 1 - \pi\left(\alpha, 1-p\right)   .
\label{Equ:pi-Symmetry}
\end{equation}
This proves that 
\[
p^2 < \pi\left(\alpha, p\right) < 2p -p^2
\]
as well as $\pi\left(\alpha, \sfrac{1}{2}\right) = \sfrac{1}{2}$.
Thus properties (ii) and (iii) of $\pi$ hold.

Moreover from the expression (\ref{Equ:pi-Expression}) it follows that for every fixed 
$\alpha \in \left(0,1\right)$ the limiting marginal $\pi$ as a function of $p$ is an
increasing limit of polynomials in $p$. This implies 
that $p \mapsto \pi\left(\alpha, p \right)$ is \emph{lower semi-continuous} \cite{RuRC87}.
But because of the identity (\ref{Equ:pi-Symmetry}) for the same reason it is also 
\emph{upper semi-continuous}. This proves that that $\pi$ as a function of $p$ is 
continuous, establishing the property (i). 

Finally, we show property (iv). For this fix $\alpha \in \left(0, 1\right)$ and
notice that from the expression (\ref{Equ:pi-Expression}), since all the summands are non-negative, we have
\begin{equation}
\frac{\pi\left(\alpha, p\right)}{p} 
\geq
\frac{1}{p} \alpha \left(1-\alpha\right)
\sum_{t=0}^{\infty} \left(\Pbb_t\left(BRBB\right) + \Pbb_t\left(BBRB\right)\right)  .
\label{Equ:pi-Derivative-1}
\end{equation}
Now fix $t \geq 0$ and consider the probability $\Pbb_t\left(BRBB\right)$. Because of translation invariance
without loss of any generality,  we may assume that the configuration we are considering is at the locations
$\left(-1,0,1,2\right)$. Now notice that because the dynamics depends only on the nearest neighbours so
$\Pbb_t\left(BRBB\right)$ depends on the initial configuration at the locations in the interval $\left[-t-1, t+2\right]$. So 
without loss of genrailty we may assume that outside the interval $\left[-t-1, t+2\right]$, at every location the colour of the chameleons are
blue ($B$). So we may write
\begin{equation}
\Pbb_t\left(BRBB\right) = p \left(1-p\right)^{2t+3} p_{11}^{(t)} \left(\alpha\right) + o\left(p^2\right)  
\label{Equ:Fundamental-Chain-Decomp}
\end{equation}
where the terms in $o\left(p^2\right)$ are all non-negative and $p_{11}^{(t)}\left(\alpha\right)$ is the sum over all locations $x \in \left[-t-1, t+2\right]$ of the
probability of obatining exactly one $R$ chameleon at location $0$ at time $t$ having started  
 at time $0$ with exactly one red ($R$) chameleon at  location $x$, 
Observe also that $\Pbb_t\left(BBRB\right)$ has exactly  the same representation as $\Pbb_t\left(BRBB\right)$.

Now let us consider the case when we start with exactly one $R$ chameleon at some location $x \in \Zbold$ and all
other chameleon of color $B$. For this let $L_t$ be the number of red chameleons at time $t$ and
$X_t$ be the position of the leftmost red chameleon at time $t$. These two quantities are well defined for
our Markov chain. Thus we get
\begin{eqnarray}
p_{11}^{(t)} 
 & = & \sum_{x=-t-1}^{t+2} \bP\left(L_t = 1, X_t = 0 \,\Big\vert\, L_0 = 1, X_0 = x\right) \nonumber \\
 & = & \sum_{x=-t-1}^{t+2} \bP\left(L_t = 1, X_t = -x \,\Big\vert\, L_0 = 1, X_0 = 0\right) \nonumber \\
 & = & \bP\left(L_t = 1 \,\Big\vert\, L_0 = 1, X_0 = 0\right)
       \label{Equ:New-Chain-Definition}
\end{eqnarray}
where the second equality follows because of the translation invariance of the measure while the last 
follows because if $X_0 = 0$ then $X_t \in \left[-t, t\right]$ with probability one. 

Now it follows easily that starting with exactly one $R$ chameleon at the origin the stochastic 
process $\left(L_t\right)_{t=0}^{\infty}$ is a Markov chain with state-space 
$\left\{0, 1, 2, \ldots \right\}$ starting at $L_0=1$ and with absorbing state $0$. The
transition matrix $P := \left(\left( p_{ij} \right)\right)$ is given by
\begin{equation}
p_{ij} 
= \left\{
  \begin{array}{ll}
  2 \alpha \left(1-\alpha\right)^2 + \alpha^2 \left(1-\alpha\right) & \mbox{if\ } i=1, j=0 \\
  1 - 3 \alpha \left(1- \alpha\right)                               & \mbox{if\ } i=1, j=1 \\
  2 \alpha^2 \left(1-\alpha\right)                                  & \mbox{if\ } i=1, j=2 \\
  \alpha \left(1 - \alpha\right)^2                                  & \mbox{if\ } i=1, j=3 \\
  \bP\left(Z_1 + Z_2 = j-i\right)                                   & \mbox{if\ } i \geq 2 \\
  0                                                                 & \mbox{otherwise}
  \end{array}
  \right. 
\label{Equ:Length-Chain-Transition-Matrix}
\end{equation}
where $Z_1, Z_2$ are i.i.d. random variables with 
$\bP\left(Z_1 = -1\right) = \bP\left(Z_1 = 1\right) = \alpha \left(1-\alpha\right)$ and 
$\bP\left(Z_1 = 0\right) = 1 - 2 \alpha \left(1-\alpha\right)$.

Now using equations (\ref{Equ:pi-Derivative-1}), (\ref{Equ:Fundamental-Chain-Decomp}) it follows that
\begin{eqnarray*}
\liminf_{p \rightarrow 0} \frac{\pi\left(\alpha, p\right)}{p} 
 & \geq & 2 \alpha \left(1-\alpha\right) 
          \liminf_{p \rightarrow 0} \sum_{t=0}^{\infty} \left(1 - p \right)^{2t+3} p_{11}^{(t)} \\
 & \geq & 2 \alpha \left(1-\alpha\right) \sum_{t=0}^{\infty} p_{11}^{(t)} \\
 & \geq & 1   .
\end{eqnarray*}
Here we note that the second inequality follows from Fatou's Lemma and, noting that 
$\sum_{t=0}^{\infty} p_{11}^{(t)} = \bE\left[ \# \mbox{ of returns to state }1 \,\Big\vert\, L_0 = 1 \right]$, 
 the last 
inequality follows from Theorem \ref{Thm:Length-Chain-Lower-Bound} of Section \ref{Sec:Technical}.

This completes  the proof of Theorem \ref{Thm:Main}. $\qed$

\section{Translation Invariant Starting Distribution}
\label{Sec:Translation-Invariant}
Exactly as in the case of the one directional model here too we observe 
from the proof of the convergence that (\ref{Equ:Convergence-to-All-Same-Color}) goes through without any 
change for any $\Pbb_0$ which is  translation invariant. Thus we get the
following generalization of Theorem \ref{Thm:Main}. 
\begin{Theorem}
\label{Thm:Main-General}
Let $\{\xi_i(0) : i \in \mathbb Z\}$ be random variables which are translation invariant and let
$\bP(\xi(0) = R) = p = 1-\bP(\xi_i(0) = B)$. For the two sided model with
$p_R = p_B = \alpha \in \left(0,1\right)$ we have
\begin{equation}
\Pbb_t \cd \bar{\pi} \delta_{\bR} + \left(1-\bar{\pi}\right) \delta_{\bB}  ,
\end{equation}
where $\bar{\pi}$ depends on the initial distribution $\Pbb_0$ and as well as $\alpha$. 
\end{Theorem}
The following corollary is also an immediate consequence  
\begin{Corollary}
\label{Cor:Characterization-Stationary-Measures}
For the two sided neighbourhood model with $p_R = p_B = \alpha \in \left(0,1\right)$ 
the only translation invariant \emph{stationary} measures are of the form
\[
\kappa \delta_{\bR} + \left(1-\kappa\right) \delta_{\bB} 
\]
for some $0 \leq \kappa \leq 1$. 
\end{Corollary}

Here it is worthwhile to mention that it is unlikely that this chain has a stationary distribution which is not
translation invariant, but we have not explored in that direction.

\section{Some Technical Results}
\label{Sec:Technical}
In this section we prove some technical results which have been used in the proofs in the previous sections. 

\begin{Theorem}
\label{Thm:General-Conv}
Let $\left(\Pbb_t\right)_{t \geq 0}$ be a sequence of translation invariant measures on
$\left\{R, B\right\}^{\Zbold}$ such that the following conditions hold
\begin{itemize}
\item[(i)]   $\lim_{t \rightarrow \infty} \Pbb_t\left(R\right)$ exists,
\item[(ii)]  $\lim_{t \rightarrow \infty} \Pbb_t\left(RR\right)$ exists, and
\item[(iii)] for all $k \geq 1$
             $\lim_{t \rightarrow \infty} \Pbb_t\left(B\bR_kB\right) 
             = 0 =
             \lim_{t \rightarrow \infty} \Pbb_t\left(R\bB_kR\right)$,
\end{itemize}
then
\begin{equation}
\Pbb_t \cd a \delta_{\bR} + \left(1-a\right) \delta_{\bB} \, \mbox{ as } t \to \infty   ,
\end{equation}
where $a := \lim_{t \rightarrow \infty} \Pbb_t\left(R\right)$.  
\end{Theorem}

\proof Let $a := \lim_{t \rightarrow \infty} \Pbb_t\left(R\right)$ and
$b := \lim_{t \rightarrow \infty} \Pbb_t\left(RR\right)$. To prove the result it is enough
to show that $a = b$. This is because then
$\Pbb_t\left(RB\right) = \Pbb_t\left(BR\right) = \Pbb_t\left(RR\right) - \Pbb_t\left(R\right) 
\rightarrow 0$ as $t \rightarrow \infty$. 

Now to show $a=b$ we first observe that 
\[
\Pbb_t\left(RB\right) = \Pbb_t\left(R\right) - \Pbb_t\left(RR\right) \rightarrow a - b  .
\]
Now for $k \geq 1$, 
\[ 
\Pbb_t\left(\bR_{k+1}B\right) = \Pbb_t\left(\bR_kB\right) - \Pbb_t\left(B\bR_kB\right)  .
\]
Thus under assumption (iii) it follows by induction that
\begin{equation}
\Pbb_t\left(\bR_kB\right) \rightarrow a - b .
\label{Equ:General-Conv-1}
\end{equation}

Also for any $k \geq 1$ we have
\[
\Pbb_t\left(\bR_{k+1}\right) = \Pbb_t\left(\bR_k\right) - \Pbb_t\left(\bR_kB\right)  .
\]
Thus it follows by induction that
\begin{equation}
\Pbb_t\left(\bR_k\right) \rightarrow \left(k-1\right) b - \left(k-2\right) a  .
\label{Equ:General-Conv-2}
\end{equation}

From equation (\ref{Equ:General-Conv-1}) and (\ref{Equ:General-Conv-2}) it follows that
\begin{equation}
b \leq a \leq \sfrac{k-1}{k-2} b \mbox{ for all } k \geq 3  .
\label{Equ:General-Conv-3}
\end{equation}
This proves that $a=b$ completing the proof. $\qed$

\begin{Theorem}
\label{Thm:Length-Chain-Lower-Bound}
Let $\left(L_t\right)_{t \geq 0}$ be a Markov chain on the state-space
$\left\{0, 1, 2, \ldots \right\}$ with transition matrix $P = \left(\left(p_{ij}\right)\right)$ as given
in (\ref{Equ:Length-Chain-Transition-Matrix}) and $L_0=1$. Then
\begin{equation}
\bE\left[\# \mbox{\ of returns to the state\ } 1 \,\Big\vert\, L_0 = 1 \right] 
> \frac{1}{2 \alpha \left(1- \alpha\right)}
\label{Equ:Lower-Bound}
\end{equation}
\end{Theorem}

\proof Let $f_{11}^{\star} := \bP\left(L_t = 1 \mbox{\ for some\ } t \geq 1 \,\Big\vert\, L_0 = 1\right)$ then
from standard Markov chain theory \cite{FellerII71} it follows that
\begin{equation}
\bE\left[\# \mbox{\ of returns to the state\ } 1 \,\Big\vert\, L_0 = 1 \right]
=
\frac{1}{1- f_{11}^{\star}}
\label{Equ:f11-Star-Representation}
\end{equation}

Moreover we can also write,
\begin{equation}
f_{11}^{\star} = p_{11} + p_{12} f_{21}^{\star} + p_{13} f_{31}^{\star}  ,
\label{Equ:f11-Star-Decomposition}
\end{equation}
where 
$f_{k1}^{\star} := 
\bP\left(L_t = 1 \mbox{\ for some\ } t \geq 1 \,\Big\vert\, L_0=k \right)$ for 
$k \in \left\{ 2, 3 \right\}$.

Now let $\left(\left(\bar{p}_{ij}\right)\right)$ be a new Markov chain on the same state-space
$\left\{0, 1, 2, \ldots \right\}$ such that both $0$ and $1$ are absorbing states and 
$\bar{p}_{ij} = p_{ij}$ for all $i \geq 2$. Let $u_k$ be the probability of getting absorbed in the
state $1$ for this new chain when started at state $k$. Then it is easy to see that
$f_{k1}^{\star} = u_k$ for any $k \geq 2$.

From definition $u_0 = 0$ and $u_1 = 1$. Moreover it is easy to see that 
\begin{equation}
u_k = \beta_2 \left(u_{k-2} + u_{k+2}\right) + \beta_1 \left(u_{k-1} + u_{k+1}\right) + \beta_0 u_k \,,
\label{Equ:u}
\end{equation}
where $\beta_2 = \theta^2$, $\beta_1 = 2 \theta \left(1-2\theta\right)$ and $\beta_0 = 1 - \beta_1 - \beta_2$
and $\theta := \alpha \left(1-\alpha\right)$. 
The characteristic polynomial of this difference equation is given by
\begin{equation}
\lambda^2 = \beta_2 \left(\lambda^4 + 1\right) + \beta_1 \left(\lambda^3 + \lambda\right) + \beta_0 \lambda^2  .
\label{Equ:u-Characteristic}
\end{equation}
It then follows easily that this has three real roots, $\gamma_0 = 1$ with multiplicity $2$ and
$-\gamma_1$ and $-\gamma_2$ such that $\gamma_2 > 1 > \gamma_1 > 0$. 
So a general solution of (\ref{Equ:u}) is given by
\begin{equation}
u_k = C_1 + C_2 k + C_3 (-\gamma_1)^k + C_4 \left(-\gamma_2\right)^k  .
\end{equation}
But our $u_k$'s are probability and hence are in $\left[0,1\right]$, so we must have $C_2 = C_4 = 0$. But 
because of the initial conditions $u_0=0$ and $u_1=1$ it follows that
\begin{equation}
u_k = \frac{1}{1+\gamma_1} - \frac{\left(-\gamma_1\right)^k}{1+\gamma_1}  .
\label{Equ:u-Expression}
\end{equation}
In particular
\begin{equation}
u_2 = 1 - \gamma_1 \mbox{\ and\ } u_3 = 1 -\gamma_1 + \gamma_1^2  .
\label{Equ:u2-u3}
\end{equation}
Going back to the characteristic equation (\ref{Equ:u-Characteristic}) we determine that
\begin{equation}
1 - \gamma_1 = h\left(\theta\right)  ,
\end{equation}
where 
\begin{equation}
h\left(\theta\right) = \frac{\sqrt{1-2 \theta} - \left(1 - 2\theta\right)}{\theta}  .
\label{Equ:h-Definition}
\end{equation}

Now to complete the proof we need to show that
\[
f_{11}^{\star} > 1 - 2 \alpha \left(1- \alpha\right)
\]
which is equivalent to showing
\[ h\left(\theta\right) \left(3 \alpha - 1\right) - \alpha 
                     + \left(1-\alpha\right)\left(h\left(\theta\right)\right)^2 > 0
\]
where $\theta = \alpha \left(1-\alpha\right)$ and $h$ is as defined in (\ref{Equ:h-Definition}).

The rest of the proof is simple calculus and some exact calculations; for completeness 
we provide the essential details. 
From the definition of $h$ one can show easily by Taylor expansion that
\begin{equation}
1 - h\left(\theta\right) \leq \sfrac{\theta}{2} + c \theta^2  ,
\label{Equ:h-bound}
\end{equation}
where $c = 2.18$.  Also it is not difficult to show that
\begin{equation}
h\left(\theta\right) \geq 2 \left(\sqrt{2} - 1\right) > 0.8  ,
\end{equation}
here we note that $\theta = \alpha \left(1-\alpha\right) \in \left[0,\sfrac{1}{4}\right]$. 

Finally, 
\begin{eqnarray*}
\lefteqn{h\left(\theta\right) 
          \left(3 \alpha - 1\right) - \alpha + \left(1-\alpha\right)\left(h\left(\theta\right)\right)^2 }\\
 &  =   & \alpha \left(2 h\left(\theta\right) - 1\right) 
          - \left(1-\alpha\right) h\left(\theta\right) \left(1-h\left(\theta\right)\right) \\
 & \geq & \alpha \left(2 h\left(\theta\right) - 1\right)       
          - \left(1-\alpha\right) h\left(\theta\right) \theta \left(\sfrac{1}{2} + c \theta\right) \\
 &  =   & \alpha \left\{ 
          h\left(\theta\right) \left( 2 - \left(1-\alpha\right)^2 \left(\sfrac{1}{2} + c \theta\right)\right)
          - 1 \right\} \\
 & \geq & \alpha \left\{ h\left(\theta\right) 
          \left( 2 - \sfrac{1}{2} - c \alpha \left(1-\alpha\right)^3 \right) - 1 \right\} \\
 & \geq & \alpha \left\{ 0.8 \left(\sfrac{3}{2} - 2.18 * \sfrac{27}{256}\right) - 1\right\} \\
 & \geq & \alpha \left(1.01 - 1\right) \\
 &  >   & 0.
 \end{eqnarray*}
Here we use the fact that $\alpha \left(1-\alpha\right)^3 \leq \sfrac{27}{256}$.
This completes the proof. $\qed$

\section{Coexistence of the Two Colours}
\label{Sec:Coexistence}
We observe that in our model even in the critical case only one of the colours survives at the limit provided we 
start with a translation invariant distribution. In this case it is also easy to see that if we start with all
blue chameleons on the negative integers and all red chameleons at the non-negative integer locations then
at the limit with probability $1/2$ we will have an all red configuration and with probability $1/2$ it will be
an all blue configuration. So once again no coexistence. This is because in this case the interface
between the blue and red chameleons will perform a discrete time symmetric random walk with holding and hence
the result will follow from standard local limit theorem. 

It is possible though to get coexistence of the two colours by making the success probability $\alpha$ inhomogeneous,
that is to depend on time. Suppose that $\alpha_t$ is the probability of success of the coin toss of any colour at time $t$ and we
start with a configuration of all blue chameleons on the negative integers while all red chameleons at the 
non-negative integers. Let $X_t$ be the position of the left-most red chameleon. Then it is easy to see that
$X_t := Y_1 + Y_2 + \cdots Y_t$ where $\left(Y_i\right)_{i \geq 1}$ are independent and $Y_i$ follows a 
distribution on $\left\{-1, 0 , 1\right\}$ with 
$\bP\left(Y_i = -1\right) = \bP\left(Y_i = 1\right) = \alpha_i \left(1- \alpha_i\right)$ and
$\bP\left(Y_i = 0\right) = 1 - 2 \alpha_i \left(1 - \alpha_i\right)$. 
So by Kolmogorov's Three Series Theorem \cite{Bill95} the sequence of random variables 
$\left(X_t\right)_{t=0}^{\infty}$ converges a.s. if and only if 
\begin{equation}
\sum_{t=1}^{\infty} \alpha_t \left(1 - \alpha_t\right) < \infty
\label{Equ:Condition-Coexistence}
\end{equation}
Thus if (\ref{Equ:Condition-Coexistence}) is satisfied then there will be both colours present 
at the limit. 

This is intuitively clear since under the condition (\ref{Equ:Condition-Coexistence}) for large enough $t$ one of
$\alpha_t$ or $\left(1-\alpha_t\right)$ is ``small'' and hence there will either be a large number of failures
or large number of successes of the coin tosses, and  in either case, no change is expected. 

It is of course more interesting to study this critical model on higher dimensions with 
homogeneous success probability and to explore the possibility of coexistence in that case. Unfortunately our method does not help in that case. In fact in higher dimension it is not even clear a condition like $p_R > p_B$ is
good enough to get all red configuration at the limit.

\section*{Acknowledgement}
The authors would like to thank Kalyan Chatterjee and Debraj Ray for useful discussions. The work of Rahul Roy
was supported by a grant from Department of Science and Technology, Government of India. 

\bibliographystyle{plain}

\bibliography{Learning-from-Neighbours}

\end{document}